\documentclass[a4paper]{article}
\usepackage{amsmath}
\usepackage{bbm}
\usepackage{mathrsfs}
\newcommand{\rf}[1][]{\textup{\eqref{#1}}}
\newcommand{\Rda}{\hat{R}}

\def\abar{{\overline{a}}}
\newcounter{beh}[section]
\def\thebeh{\arabic{beh}}

\newenvironment{thm}{\noindent\refstepcounter{beh}{\bf Theorem\
\thebeh.} \em }{\em}

\newenvironment{prp}{\noindent\refstepcounter{beh}{\bf Proposition\
\thebeh.} \em }{\em}

\newenvironment{Proof}{\noindent{\it Proof}.}{\ \hspace*{\fill}
\rule{1ex}{1ex}}

\def\id{\mathrm{id}}
\def\End{\mathrm{End}}

\def\dd{\mathrm{d}}
\def\ot{\otimes}

\def\AAA{\mathscr{A}}
\def\SS{\mathscr{S}}
\def\LL{\mathscr{L}}
\def\ant{A}
\def\II{{\scriptscriptstyle{\mathrm{I}{\hspace{-1pt}}\mathrm{I}}}}
\def\I{\scriptscriptstyle{\mathrm{I}}}

\def\M{{M}}
\def\S{\mathcal{S}}
\def\AA{\mathcal{A}}
\def\B{\mathcal{B}}
\def\CC{\mathcal{C}}

\def\C{\mathbbm{C}}
\def\N{\mathbbm{N}}
\def\NO{\mathbbm{N}_0}
\def\R{\mathbbm{R}}

\def\slqn{SL_q(N)}
\def\glqn{GL_q(N)}

\def\nn{\nonumber}

\def\tsum{\textstyle \sum}

\newcommand{\qm}{q^{-1}}

\newcommand{\lam}{{\lambda}}
\newcommand{\QM}{Q}
\newcommand{\QP}{Q_+}

\def\ve{\varepsilon}

\def\vr{\varrho}
\def\vN{\varUpsilon}
\def\vO{\varOmega}
\def\om{\omega}

\def\rac{\triangleleft}

\newcommand \fettc{{\scriptstyle{\tt c}}}

\def\RR{\mathcal{R}}
\def\AAi{\AA_{\mathrm{inv}}}
\def\RRi{\RR_{\mathrm{inv}}}
\def\ii{\mathrm{i}}
\def\ll{\ell}
\def\rr{\mathrm{r}}

\def\rrr{r}
\def\rp{\rrr_+}
\def\rmi{\rrr_-}
\def\tttt{\theta_\tau}

\def\adr{{\mathrm{Ad}_\rr}}

\def\dl{{\Delta_\ll}}
\def\dr{{\Delta_\rr}}

\def\dlx{\dl\colon\Gamm\to\AA\ot\Gamm}
\def\drx{\dr\colon\Gamm\to\Gamm\ot\AA}

\newcommand{\dc}{(\Gamm,\dd)}
\newcommand{\Gamm}{\varGamma}

\def\gpmz{\Gamm_{\pm,z}}

\def\gpz{\Gamm_{+,z}}
\def\gmz{\Gamm_{-,z}}

\newcommand\gl{{\Gamm_\ll}}
\newcommand\gr{{\Gamm_\rr}}
\newcommand\gi{{\Gamm_\ii}}

\def\gon{{\Gamm^n}}

\newcommand\ota{{\otimes}_{\!{\scriptscriptstyle\mathcal{A}}}}

\def\sig{\sigma}
\newcommand\gten{{\Gamm^{\otimes}}}
\newcommand\gtenn{{\Gamm^{\otimes n}}}
\newcommand\gtenk{{\Gamm^{\otimes k}}}

\newcommand\ww{{\mbox{$\scriptscriptstyle{W}$}}}
\newcommand\gd{{\Gamm^{\land}}}
\newcommand\gdw{{_\ww{\Gamm^{\land}}}}
\newcommand\gds{{_s{\Gamm^{\land}}}}
\newcommand\gdu{{_u{\Gamm^{\land}}}}

\def\ul{\underline}
\def\JJ{{_\ww J}}
\def\Js{{_s J}}
\def\Ju{{_u J}}

\newcommand\gdwl{{_\ww{\Gamm^{\land}_\ll}}}

\newcommand\gdsl{{_s{\Gamm^{\land}_\ll}}}
\newcommand\gdsi{{_s{\Gamm^{\land}_\ii}}}
\newcommand\gdul{{_u{\Gamm^{\land}_\ll}}}
\newcommand\gdwi{{_\ww{\Gamm^{\land}_\ii}}}

\newcommand\gdwik{{_\ww{\Gamm^{\land k}_\ii}}}
\newcommand\gdwin{{_\ww{\Gamm^{\land n}_\ii}}}

\begin{document}

\author{
I.~Heckenberger\thanks{e-mail: heckenbe@mathematik.uni-leipzig.de}
\,\,and
 A.~Sch{\" u}ler\thanks{Supported by the Deutsche Forschungsgemeinschaft.
schueler@mathematik.uni-leipzig.de}
\\
{\small\sl Institute of Mathematics, University of Leipzig}
\\
{\small\sl Augustusplatz 10, 04109 Leipzig, Germany}
}   

\title{Higher Order Differential Calculus on $SL_q(N)$}

\date{}
\maketitle
\begin{abstract}

Let $\Gamm$ be a bicovariant first order differential calculus on
a Hopf algebra $\AA$. There are three possibilities to construct
a differential $\NO$-graded Hopf algebra $\gd$ which contains $\Gamm$ as its
first order part.
In all cases $\gd$ is a quotient  $\gd=\gten /J$ of the tensor algebra by some suitable ideal.
We distinguish	three possible choices $\Ju$, $\Js$, and $\JJ$, where the first one generates
the universal differential calculus (over $\Gamm$) and the last one is
Woronowicz' external algebra.
\\
Let $q$ be a transcendental complex number and let $\Gamm$ be one
of the $N^2$-dimensional bicovariant first order differential calculi on the
quantum group $\slqn$. Then for $N\ge3$ the three ideals coincide.
For Woronowicz' external algebra we calculate the dimensions
of the spaces of left-invariant and bi-invariant $k$-forms. In this case each
bi-invariant form is  closed.
\\
In case of $4D_\pm$ calculi on $SL_q(2)$ the universal calculus is
strictly larger
than the other two calculi. In particular, the bi-invariant 1-form is {\em not} closed.
\vspace{1cm}\noindent
\end{abstract}
{\bf 1.} It was S.\,L.\,Woronowicz who provided a general framework for bicovariant differential
calculus on quantum groups  \cite{Wo2}.
Covariant first order differential calculi (abbreviated FODC) were
constructed, studied and classified by many authors, see for instance
\cite{Wo1,CSWW,J,SS1}.
Despite the rather extensive
literature on bicovariant FODC the corresponding differential calculi of higher
order forms have been treated only in few cases. In \cite{Brz}
it was pointed out that the exterior algebra
$\gds$ is a differential graded Hopf algebra.
In \cite{gr} the de Rham cohomology of $4D_\pm$ calculi $\gds$ was calculated.
After giving the talk at the conference we were informed about a paper
of P.\,P.\,Pyatov and L.\,D.\,Faddeev \cite{PF} where the algebras
$\gdsl$ and $\gdsi$ on $\glqn$ were also investigated.
\\
The purpose of this paper is to compare three  possible constructions
of exterior algebras to a given bicovariant FODC. Consider the tensor
algebra $\gten$. Let $\Ju$ be  the two-sided ideal generated by
the elements $\sum_{(r)}\omega(r_{(1)})
\otimes\omega(r_{(2)})$, $r\in \RR$, where $\omega(a)=
\sum_{(a)}S(a_{(1)}){\rm d}(a_{(2)})$ and
$\RR=\ker\varepsilon\cap\ker\omega$.
Let $\Js$ be generated by $\ker(I-\sig)$, where $\sig$ is the braiding
of $\Gamm\ota\Gamm$.
Finally let  $\JJ=\bigoplus_{k\ge 2}\ker \ant_k$, where $\ant_k$
is the $k$th antisymmetrizer
constructed {}from the braiding $\sig$.
Define the exterior algebras $\gdu=\gten/\Ju$, $\gds=\gten/\Js$,
 and $\gdw=\gten/\JJ$
and call them {\em universal exterior algebra},
{\em ``second antisymmetrizer" exterior algebra} and {\em Woronowicz'
external algebra}, respectively.
\\
The first main result stated in Theorem\,\ref{left} and Theorem\,\ref{bi} is
the calculation of Poincar\' e polynomials  of the subalgebras of
left-invariant
and bi-invariant forms. Both coincide with the Poincar\' e
polynomials of left-invariant
and bi-invariant forms on the classical Lie group $GL(N,\R)$, respectively. Details
of the proofs
will be found in \cite{Sch}.
The second main result stated in Theorem\,\ref{usw}
compares the three algebras of left-invariant forms. Suppose $q$  to be
transcendental,
$N\ge 3$, $\AA=\slqn$ and $\Gamm=\Gamm_{\pm,k}$. Then $\gdul$,  $\gdsl$,
and $\gdwl$ are isomorphic differential graded algebras.
For $\AA=SL_q(2)$ and $\Gamm=4D_\pm$ however
$\gdul$ is strictly larger than $\gdsl$.
\\
In Section\,2 we recall the basic
definitions and preliminary facts on bicovariant bimodules,  bicovariant
 FODC and the  construction of
bicovariant FODC on $\glqn$ and  $\slqn$.
In the last part we state the main results and prove
Theorem\,\ref{usw}.
\\[.5cm]
{\bf 2.}
Let $\AA$ be a Hopf algebra.
We denote by $\Delta$ the comultiplication,
$\ve$ the counit, and $S$ the antipode.  The linear span of a set $\{a_i:i\in I\}$
is denoted by $\langle a_i: i\in I\rangle$. As usual $v^\fettc$ ($f^\fettc$)
denotes the contragredient corepresentation (representation) of $v$ ($f$).
We use Sweedler's notation for the coproduct $\Delta(a)=\sum a_{(1)}
\ot a_{(2)}$,
for left coactions  $\dl(v)=\sum v_{(-1)}\ot v_{(0)}$, and for right
coactions  $\dr(v)=\sum v_{(0)}\ot v_{(1)}$.

\noindent{\em Bicovariant bimodules.}
Let $(\Gamm,\dl,\dr)$  be a bicovariant bimodule  \cite{Wo2} over $\AA$.
 An element
$\omega\in\Gamm$ is called {\em left-invariant} (resp. {\em right-invariant})
if $\dl(\omega)=1\ot\omega$ (resp. $\dr(\omega)=\omega\ot 1$).
The linear space of left-invariant (resp. right-invariant) elements
of $\Gamm$ is denoted by $\gl$
(resp. $\gr$). The elements of $\gi=\gl\cap\gr$ are called {\em bi-invariant}.
The structure of bicovariant bimodules over $\AA$ is completely characterized
by Theorems 2.3 and 2.4 in \cite{Wo2}. We recall the corresponding result:
\\
Let $(\Gamm,\dl,\dr)$ be a bicovariant bimodule over $\AA$ and let
$\{\omega_i\}$ be a linear basis of $\gl$.
Then there exist matrices $v=(v^i_j)$
and $f=(f^i_j)$ of elements $v^i_j\in \AA$ and of functionals $f^i_j$ on $\AA$
such that:
\\
(i) $\sum Sa_{(1)}\omega_i a_{(2)}
=\sum_n f^i_n(a)\omega_n$, $a\in \AA$ and $\dr(\omega_i)=
\sum_n\omega_n\ot v^n_i$.
\\
(ii) $v$ is a corepresentation and $f$ is a representation of $\AA$.
\\
(iii) $\sum_r v^r_i(a\ast f^r_j)=\sum_r(f^i_r\ast a)v^j_r$, $a\in\AA$.\\
The set $\{\omega_i\}$ is a free left module basis of $\Gamm$.
 We have set $a\ast f=
\sum a_{(1)}f(a_{(2)})$ and $f\ast a=\sum a_{(1)}f(a_{(2)})$.
Conversely, if $\{\omega_i\}$
is a basis of a certain vector space $\Gamm_0$ and if $v$ and $f$
are matrices satisfying (ii)
and (iii) then there exists a unique bicovariant bimodule $\Gamm$ such
that $\gl=\Gamm_0$ and (i) holds. In this situation we simply
write $\Gamm=(v,f)$.
\\
Let $\Gamm_1=(v_1,f_1)$ and $\Gamm_2=(v_2,f_2)$ be bicovariant bimodules.
The pair $(v_1\ot v_2, f_1\ot f_2)$ is also a bicovariant bimodule
which will be denoted
by $\Gamm_1\ota \Gamm_2$.
In this way
$\gtenn=\Gamm\ota\cdots\ota\Gamm$ ($n$ factors), $n\ge2$, and
$\gten=\bigoplus_{n\ge 0}\gtenn$, where
$\Gamm^{\ot 0}=\AA$, $\Gamm^{\ot 1}=\Gamm$,
become bicovariant bimodules as well.\\
For a bicovariant bimodule $\Gamm$
the space $\Gamm_\ll$ becomes a right $\AA$-module via
 $ \rho\rac a=\sum Sa_{(1)}
\rho a_{(2)}$. Moreover,
$(\rho\,\omega)\,\rac a=\sum_{(a)}(\rho\rac a_{(1)})(\omega\rac a_{(2)})$
in both $\gten$ and $\gd$.
\\
\noindent{\em  Shuffle decomposition and
lift into braids.} Let $s_1,\dots,s_{k-1}$
 denote the nearest neighbor transpositions
 of the symmetric group $\S_k$. The elements of $\CC_{ki}=
\{p\in\S_k\colon p(m)<p(n)\text{ for } 1\le m<n\le i \text{ and }
i+1\le m<n\le k\}$ are called {\em shuffle} permutations. Each $p\in\S_k$
admits a unique representation $p=p_1p_2p_3$ where $p_1\in\CC_{ki}$
and $p_2,\,p_3\in \S_k$ leave the numbers $\{i+1,\dots,k\}$
respectively $\{1,\dots ,i\}$ pointwise fixed. Moreover, $\ell(p)=
\ell(p_1)+\ell(p_2)+\ell(p_3)$, where $\ell$ denotes the  length-function
on $\S_k$. Let $b_1,\dots,b_{k-1}$ be the generators of the
braid group $\B_k$.
For a reduced expression $w=s_{i_1}\cdots s_{i_r}$, i.\,e. $\ell(w)=r$,
define the   braid $b_w=b_{i_1}\cdots b_{i_r}$.
Obviously,  $b_w$ is independent of the choice of the reduced
expression $w$. We have
$b_{vw}=b_vb_w$ for $\ell(vw)=\ell(v)+\ell(w).$
In particular, this equation  applies to the shuffle decomposition
$b_w=b_{p_1}b_{p_2}b_{p_3}$. Define the antisymmetrizer $\AAA_k$ and
shuffle sums $\AAA_{ki}$ in the group algebra $\C(\B_k)$ as follows
\begin{equation}
\label{symdefi}
\AAA_k=\tsum_{w\in \S_k}(-1)^{\ell(w)} b_w, \qquad
\AAA_{ki}=\tsum_{w\in \CC_{ki}}(-1)^{\ell(w)}b_w.
\end{equation}
By the shuffle decomposition we obtain for $i<k$
\begin{equation}
\AAA_k=\AAA_{ki}(\AAA_i\ot \AAA_{k-i}).    \label{sym}
\end{equation}
Now we recall the construction of the
external algebra $\gdw$ due to \cite{Wo2}.
There is a unique isomorphism $\sig$
 of bicovariant bimodules
$\sig\colon\Gamm\ota\Gamm\to\Gamm\ota\Gamm$ with
$\sig(\omega\ota\eta)=\eta\ota\omega$,
$\omega\in\Gamm_{\ll}$, $\eta\in\Gamm_{\rr}$ called the {\em braiding}.
Since $\sig$ fulfills the braid equation, it  can be
extended to a representation
$\sig$ of the group algebra $\C(\B_k)$ in
$\End_{{\scriptscriptstyle\mathcal{A}}}
(\Gamm^{\ot k})$.
The images $\sig(\AAA_k)$ and $\sig(\AAA_{ki})$ of the antisymmetrizer
and the shuffle sums will be denoted by
$\ant_k$ and $\ant_{ki}$, respectively.
By \rf[sym],
$
\JJ=\bigoplus_{k\ge 2}\ker \ant_k
$
is a twosided ideal in $\gten$. Since $\ant_k$ is a homomorphism of
bicovariant bimodules $\JJ$ is a bicovariant subbimodule of $\gten$.
Consequently, $\gdw=\gten/\JJ$ is a $\NO$-graded algebra and a
bicovariant bimodule over $\AA$.
\\
Recall that an $\NO$-graded algebra $H=\bigoplus_{n\ge 0} H^n$ is
called {\em $\NO$-graded Hopf algebra}
if the product in $H\ot H$ is given by
$
(a\ot b)(c\ot d)=(-1)^{ij}ac\ot bd,\quad b\in H^i,\,c\in H^j,\,a,\,d\in H,
$
and there are linear mappings $\Delta$, $\ve$, and $S$ of degree $0$
called coproduct, counit and antipode, respectively,   such that
the usual Hopf algebra axioms are fulfilled.
We need the following results, see for instance  \cite{KS},
Sections\,13.2 and 14.4.
 Let $\Gamm$ be  a bicovariant bimodule over $\AA$. The Hopf algebra
 structure of $\AA$
uniquely extends to an $\NO$-graded Hopf algebra structure
on $\gten$ such that
for $\omega\in\Gamm$:
$\Delta(\omega)=\dl(\omega)+\dr(\omega)$, $\ve(\omega)=0$, and
$S(\omega)=-\sum S(\omega_{(-1)})\omega_{(0)}S(\omega_{(1)})$.
The antipode is a graded antihomomorphism i.\,e. $S(\rho_1\ota\rho_2)
=(-1)^{kn}S\rho_2\ota S\rho_1$ for $\rho_1\in\gtenk$, $\rho_2\in\gtenn$.
Moreover  $\Ju,\Js$,  and $\JJ$ are
Hopf ideals in	$\gten$ (see \cite{Brz, KS}).
Hence, $\gdu$, $\gds$, and $\gdw$ are $\NO$-graded Hopf algebras.
Let $\B$ be an $\NO$-graded algebra, $\B=\bigoplus_{k\ge 0}\B_k$.
The formal power series $P(\B,t)=\sum_{k\ge 0}(\dim \B_k)\,t^k$,
is called the  {\em Poincar\'e series} of $\B$.

\noindent {\em Bicovariant Differential Calculus.}
A {\em first order differential calculus} over $\AA$ abbreviated
FODC  is a pair
$\dc$ of an $\AA$-bimodule $\Gamm$ and a linear mapping $\dd\colon\AA\to\Gamm$
such that $\dd(ab)=\dd a{\cdot} b+a{\cdot} \dd b$ for $a,\,b\in\AA$
and $\Gamm=\AA{\cdot}\dd\AA$.
A {\em differential graded algebra} over $\AA$ is a pair $(\Gamm,\dd)$
of an $\NO$-graded algebra $\Gamm=\bigoplus_{n\ge0}\gon$,
$\Gamm^{0}=\AA$,
and a linear mapping
$\dd\colon\Gamm\to\Gamm$ of degree $1$ such that $\dd^2=0$ and
\begin{equation}
\label{Leib}\nn
\dd(\omega\vr)=\dd\omega{\cdot}\vr +(-1)^n \omega{\cdot}\dd
\vr,\quad\omega\in\gon,\,\vr\in\Gamm.
\end{equation}
$(\Gamm,\dd)$ is called  {\em differential calculus} over $\AA$,
if $\gon=\AA{\cdot}\dd\AA\cdots\dd\AA$ ($n$ factors).
\\
A differential calculus (resp. FODC) $(\Gamm,\dd)$ is
called {\em bicovariant}
if there exist linear mappings $\dlx$ and $\drx$ such that
$(\Gamm,\dl,\dr)$ is a bicovariant bimodule and
\begin{equation}\nn
\begin{split}
\dl(\dd \omega)&=(\id\ot\dd)\dl(\omega), \\
\dr(\dd\omega)&=(\dd\ot\id)\dr(\omega),
\end{split}\, \text{for }\, \omega\in\Gamm,
\quad\left(\text{resp. }
\begin{split}
\dl(\dd a)&=(\id\ot\dd)\Delta(a),\\
\dr(\dd a)&=(\dd\ot\id)\Delta(a),
\end{split}\,
\text{for }\,a\in\AA.\right)
\end{equation}
{\em Adjoint corepresentation.} For $a\in \AA$ set $\adr a=
\sum a_{(2)}\ot Sa_{(1)} a_{(3)}$. The mapping $\adr\colon\AA\to\AA\ot\AA$
is a right comodule map called the {\em right adjoint corepresentation}
of $\AA$.
A crucial role play the two mappings $\omega\colon\AA\to\gl$,
$\om(a)=\sum Sa_{(1)}\dd a_{(2)}$, and $\SS\colon\AA\to\gl\ot\gl$,
$\SS(a)= \sum \om(a_{(1)})\ot\om(a_{(2)})$.
\\
According to Theorems\,1.5 and 1.8 in \cite{Wo2}  there is a one-to-one
correspondence between bicovariant FODC $\dc$ and $\adr$-invariant right
ideals $\RR$  in $\ker\ve$ given by $\RR=\ker\ve\cap\ker\om$.
The mappings $\om$ and $\SS$ are coupled by the Maurer-Cartan equation
(see  \cite{KS}, Sect.\,14, Prop.\,12):
Let $\Gamm$ be a bicovariant FODC over $\AA$ and
let $\gd=\gten/J$ be an arbitrary differential calculus over $\AA$ containing $\Gamm$
as its first order part. Then
\begin{equation}
\label{MC}
\dd(\om(a))=-\sum\om(a_{(1)})\land\om(a_{(2)})
\end{equation}
for $a\in\AA$. In particular, inserting
$a\in\RR$ into \rf[MC], the left hand side vanishes
by definition of $\RR$. Hence, $\sum \om(a_{(1)})\land\om(a_{(2)})=0$
 for $a\in\RR$
and consequently $\Ju\subseteq J$. This proves the universal property of $\gdu$.
\\
Let $(\Gamm_i,\dd_i)$, $i=1,2$, be differential graded algebras over $\AA$.
Then $(\Gamm_1\ot\Gamm_2,\dd_\ot)$ becomes a differential graded algebra
over $\AA$ if the product  in $\Gamm_1\ot\Gamm_2$ is defined
by
$
(\omega_1\ot\omega_2)(\vr_1\ot\vr_2)=(-1)^{ij}\omega_1\vr_1\ot\omega_2\vr_2,
\quad\omega_2\in\Gamm_2^{i},\,\vr_1\in\Gamm_1^{j},\,
\omega_1\in\Gamm_1,\,\vr_2\in\Gamm_2
$
and the differential $\dd_\ot$ is given by
\begin{equation}
\label{diff}  \nn
\dd_\ot(\omega_1\ot\omega_2)=\dd_1\omega_1\ot\omega_2+
(-1)^{i}\omega_1\ot\dd_2\omega_2,
\quad\omega_1\in\Gamm_1^{i},\,\omega_2\in\Gamm_2.
\end{equation}
A {\em bicovariant differential graded Hopf algebra}
 is a pair $\dc$ which is both
an $\NO$-graded Hopf algebra and a bicovariant differential graded algebra
such that the coproduct
commutes with the differentiation, i.\,e.
$\dd_\ot(\Delta(\omega))=\Delta(\dd\omega)$, for
$\omega\in\Gamm$.
Our main objects  $\gdu,\, \gds$, and $\gdw$ are  bicovariant
differential graded Hopf algebras, see \cite{KS} for details.

\noindent{\em Bicovariant FODC on Quantized Simple Lie Groups}.
We follow the method of \cite{J}  to construct
bicovariant FODC on quantizations of simple Lie groups.
Let $\AA$ be one of the quantum groups $\slqn$ or $\glqn$
as defined in \cite{FRT},
Definition\,3 and Remark\,4. Let $u=(u^i_j)_{i,\,j=1,\dots,N}$ be
the fundamental matrix corepresentation.
The deformation parameter
$q$ is assumed to be a transcendental complex number. We set $\QM=q-\qm$
and $\QP=q+\qm$.
We recall the definition of $\ell^\pm$ representations of $\AA$.
For  $x\in\C^\times$ let
 $\ell^\pm_x=((\ell_x^\pm)^i_j)$ be the $N\times N$
matrix of linear functionals
$(\ell_x^\pm)^i_j$ on $\AA=\glqn$ determined
by the properties ${\ell^{\pm}_{x}}^i_j(u^m_n)=x^{\mp 1}
(\Rda^{\pm 1})_{nj}^{im}$ and
$\ell^\pm_x\colon\AA\to\M_N(\C)$ is a unital algebra homomorphism.
For $\AA=\slqn$ we have to  assume $x^N=q$ (see also
\cite{SS1}).
We define the bicovariant
FODC $\gpmz$:
\begin{equation} \nn         \label{gpmz}
\gpz=(u^\fettc\ot u,\ell^+_x\ot\ell^{-,\fettc}_y),   \quad xy=z^{-1},\quad
\gmz=(u^\fettc\ot u, \ell^-_x\ot \ell^{+,\fettc}_y),\quad xy=z.
\end{equation}
The structure of $\gpmz$ can easily be described as follows. There exists
a basis $\{\nu^i_j:i,\,j=1,\dots,N\}$ of $(\gpmz)_\ell$ such that
the right action
and the right coaction  are given by
\begin{equation}
\label{oma}
\nu^i_j a=({\ell^\pm}^i_mS{\ell^\mp}^n_j\ast a)\nu^m_n,\quad a\in \AA
\end{equation}
and
$\dr \nu^i_j=\nu^m_n\ot(u^\fettc)^m_iu^n_j$, $i,\,j=1,\dots,N.$
The element $\nu_1=\sum_iq^{-2i}\nu^i_i$ is the unique up to scalars
bi-invariant element. Defining
\begin{equation}
\label{inner}
\dd a=\nu_1 a-a\nu_1
\end{equation}
for $a\in \AA$, the pair $(\gpmz,\dd)$
 becomes a bicovariant FODC over $\AA$.
 Transformation formulas between $\nu^i_j$ and
  Maurer-Cartan forms $\om^i_j=\om(u^i_j)$ are given
as follows.
Let $r=\rp=\qm\QM$ and $\rmi=-q^{-2N-1}\QM$.
For $n\in\N$ set $(n)_q=q^{-2}+q^{-4}+\cdots +q^{-2n}$.
Suppose $\tttt=z\bigl((N)_q+\rrr_\tau\bigr)-(N)_q$ is nonzero.
 Using matrix notation and leg numbering with
roman numbers for the objects
$\vO=(\om^i_j)$, $\vN=(\nu^i_j)$, and $u$, where
$\vN_{\I}=\vN\ot I$, $\vN_{\II}=I\ot \vN$ and $u_{\II}=I\ot u$,
we have
\begin{equation}
\label{zero}
\vO    =zr_\tau\vN +(z-1)\nu_1I\quad
\text{and}\quad\nu_1=\tttt^{-1}\om_1.
\end{equation}
For $a=(u^n_m)$ the commutation rules \rf[oma] can be written as
\begin{equation}      \label{rac}
\vN_{\I}\rac u_{\II}=z \Rda^\tau\vN_{\II}\Rda^\tau\quad
\text{and}\quad
\nu_1\rac u	=z \rrr_\tau\vN +z\nu_1 I.
\end{equation}
Now we formulate our main results.
\medskip\\
\begin{thm}     \label{left}
Let $q\in\C$ be transcendental and $\Gamm$
one of the $N^2$ dimensional bicovariant first order differential calculi
$\Gamm_{\tau,z}$, $z\in\C^\times$, $\tau\in\{+,-\}$, on
 $GL_q(N)$ or $SL_q(N)$, $z^N=q^2$. Let $\gdw$ be Woronowicz' external algebra
and let $\gdwl$ be their subalgebra of left-invariant forms.
\\
The structure of $\gdwl$ as differential graded algebra depends
only on $\tau$ but neither on $z$  nor	on  the choice of
the quantum group $GL_q(N)$ or $SL_q(N)$. The
Poincar\' e series of $\gdwl$ is
\begin{equation}  \label{poil}
P(\gdwl,t)=(1+t)^{N^2},
\end{equation}
In particular, there is a unique up to scalars least
left invariant form \mbox{$v$ of degree $N^2$.}
\end{thm}
\smallskip\\
%
\noindent\begin{thm}\label{bi}
Let $\AA$, $q$ and $\Gamm$ be as in Theorem\,\ref{left}. Let
$\gdwi$ be the subalgebra of $\gdw$ consisting of all bi-invariant
forms.
\begin{description}
\item[(i)] $P(\gdwi,t)=(1+t)(1+t^3)\cdots(1+t^{2N-1})$.
\item[(ii)]
$
\omega_k\land\omega_n=(-1)^{kn}\,\omega_n\land\omega_k,\quad
$
for\, $\omega_k\in {\gdwik}$ and $\omega_n\in{\gdwin}$.
\item[(iii)] For $\om\in\gdwi$ we have $\quad\dd \omega=0.$
\end{description}
Different bi-invariant forms represent different de Rham cohomology
classes.
\end{thm}
\medskip\\
\noindent Remark. The dimension $\dim(\gdwik)$ is equal to the number of
partitions of $k$ into a sum of pairwise different
odd numbers less than $2N$. It is also equal to the
number of  symmetric Young diagrams ($\lam=\lam'$)
with  $\lambda_1\le N$.
The least \mbox{left-invariant form $v$ is bi-invariant.}
\bigskip\\
\begin{thm}  \label{usw}
Let $\AA=\slqn$ and $q$ and $\Gamm$ as in Theorem\,\ref{left}.
 Let $\gdu$ and $\gds$ denote the
universal differential calculus and the ``second antisymmetrizer'' differential
calculus over $\Gamm$, respectively. $\gdul$ and $\gdsl$ denote their left
invariant subalgebras.
\begin{description}
\item[(i)]
$\gdsl  \cong\gdwl,\quad N\ge 2.$
\item[(ii)]
$  \gdul\cong	 \gdsl,\quad N\ge 3.$
\item[(iii)]
For $N=2$ and $\Gamm=4D_\pm$ we have
$
\quad  \gdsl\cong\gdul/\langle\omega_1\land\omega_1\rangle.
$
\end{description}
\end{thm}
\noindent Remarks. 1. For $\AA=GL_q(n|m)$ and $\Gamm=\Gamm_{\pm,1}$
(ii)  was proved in \cite{LS}.
\\
2. We expect the  de Rham cohomology of $\gdul$ in the case
$\Gamm=4D_\pm$ to be  acyclic, i.\,e. the cohomologies ${_u H}^k_\ll$
vanish for $k\ge1$.
\smallskip\\
\noindent
\begin{Proof}
(i) We recall
a result of Tsygan \cite{Tsy} on the algebra of left-invariant forms
on $\glqn$:
The algebra $\gdsl$ for $\Gamm=\Gamm_{\tau,k}$
is isomorphic to the quadratic unital complex algebra generated
by $N^2$ variables $T=(T^i_j)$ subject to the relations
\begin{equation} \label{plus}
\Rda^\tau T_\II\Rda^\tau T_\II\Rda^\tau+T_\II\Rda^\tau T_\II  =0,
\end{equation}
with  $T_\II=I\ot T$. Its Manin dual \cite{Man}
is the {\em reflection equation algebra} $\LL$
with $N^2$ generators $L=(L^i_j)$ and generating relations
\begin{equation}     \label{+}
\hat R^\tau L_{\I}\hat R^\tau L_{\I}=L_{\I}\hat R^\tau L_{\I}\hat R^\tau.
\end{equation}
Using transmutation theory invented by S.\,Majid \cite{Majid},
$\LL$ and the quantum matric bialgebra $\M_q(N)$ are isomorphic
as {\em graded linear spaces}.
Since moreover $\M_q(N)$ is a flat deformation of the algebra of polynomials
with $N^2$ indeterminates one gets, $\LL$ is also a flat
deformation  i.\,e. $P(\LL,t)=(1-t)^{-N^2}$.
Since $\LL$ provides a Poincar\' e-Birkhoff-Witt basis, $\LL$ is
Koszul by Theorem\,5.3 in \cite{Priddy}. Using Prop.\,7,\,Sect.\,9 of
\cite{Man} the Poincar\' e series of the dual quadratic
algebra  is  $P(\gdsl,t)=(1+t)^{N^2}$.
\\
Note that $\Js \subseteq\JJ$ by \rf[sym]. Since
$P(\gdwl,t)=P(\gdsl,t)$ by \rf[poil],
the ideals must be equal and we proved (i). The proofs of (ii)
and (iii) are direct consequences of Proposition \ref{su} below.
\end{Proof}
\bigskip\\
\begin{prp}    \label{su}
Let $\dc$ be one of the $N^2$ dimensional bicovariant FODC
$\Gamm_{\pm,k}$ on $\AA=\slqn$. Let $\{\nu^i_j\}$ be the
free left module basis of $\Gamm_{\pm,k}$.
\begin{description}
\item[(i)]  $\quad \nu_1\ota\nu_1\in\SS(\RR)\,\Longleftrightarrow\quad
 \gdul\cong\gdsl\quad$ {as differential graded algebras.}
\item[(ii)] $\quad\nu_1\ota\nu_1\in \SS(\RR)\quad$ for	$N\ge 3$.
\item[(iii)] $\SS(\RR)=\ker\ant_2\oplus \langle\nu_1\ota\nu_1\rangle\quad$
 for $N=2$.
\end{description}
\end{prp}
\noindent
\begin{Proof}
(i)$\Leftarrow$ is  trivial since $(I-\sig)(\nu_1\ota\nu_1)=0$
because $\nu_1$ is bi-invariant and  $\Ju=\Js$ by assumption.
\\
(i)$\Rightarrow$
Since $\Ju\subseteq\Js$ by universality, $\gdsl$ is a quotient of $\gdul$.
We show: If  $\nu_1\land\nu_1=0$ in $\gdul$, then
 the generating elements $\nu^i_j$ of $\gdul$ satisfy
\rf[plus].
For the remainder of the proof we skip the symbol $\land$
remembering that we  multiply in $\gdul$.
Our first step is to show
\begin{equation}\label{x}
\rrr_\tau\vN\vN+\vN\nu_1+\nu_1\vN=0,
\end{equation}
or more detailed $\rrr_\tau \tsum_x\nu^i_x\nu^x_j+\nu_1\nu^i_j+\nu^i_j\nu_1=0$
for  $i,\,j=1,\dots, N$.
 Using $\nu_1\nu_1=0$, \rf[inner] and the fact that each 1-form $\rho$ can be written as
$\rho=\sum a_i\dd b_i$,  one checks
\begin{equation}\label{inner1}
\dd \rho=\nu_1\rho +\rho\nu_1.
\end{equation}
Applying  \rf[MC]  to $a=u^i_j$  we obtain $\dd\om^i_j=-\sum\om^i_x\om^x_j$.
Inserting $\om^i_j$ into \rf[inner1] and comparing both equations we get in matrix notation
\begin{equation}\label{ddd}
\nu_1\vO +\vO\nu_1=\dd \vO=-\vO\vO
\end{equation}
Inserting \rf[zero] into \rf[ddd] and using $\nu_1\nu_1=0$ twice we get
\begin{equation}			       \nn
zr_\tau(\nu_1\vN+\vN\nu_1)=-z^2r_\tau^2\vN\vN-zr_\tau(z-1)(\vN\nu_1+\nu_1\vN)
\end{equation}
Subtracting the whole left hand side  we obtain \rf[x].
We complete the proof for the sample  $\tau=+$.
Applying $\rac u$ to \rf[x], using \rf[rac], $\Rda^2=\QM\Rda+I$,
and again \rf[x] (to the underlined terms) we have
\begin{align}
0=&(r\vN_{\I}\vN_{\I}+\nu_1\vN_{\I}+\vN_{\I}\nu_1)\rac u_{\II}=
(r \vN_{\I}+\nu_1I_{\I})\rac u_{\II}(\vN_{\I}\rac u_{\II}) +(\vN_{\I}\rac u_{\II})
(\nu_1I_{\I}\rac u_{\II})
 \nn\\
=&(rz\Rda\vN_{\II}\Rda+rz\vN_{\II}+z\nu_1)z\Rda\vN_{\II}\Rda
+rz\Rda\vN_{\II}\Rda\vN_{\II}+z^2\Rda\vN_{\II}\nu_1\Rda
\nn\\
=& rz^2\Rda\vN_{\II}(\QM\Rda +\ul{I})\vN_{\II}\Rda+rz^2\vN_{\II}\Rda\vN_{\II}\Rda+
\ul{z^2\Rda(\nu_1\vN_{\II}+\vN_{\II}\nu_1)\Rda}+rz^2\Rda\vN_{\II}\Rda\vN_{\II}
\nn\\
=&rz^2(\QM\Rda\vN_{\II}\Rda\vN_{\II}\Rda+\vN_{\II}\Rda\vN_{\II}\Rda+
\Rda\vN_{\II}\Rda\vN_{\II})=
rz^2\Rda (\Rda\vN_{\II}\Rda\vN_{\II}\Rda+\vN_{\II}\Rda\vN_{\II}).\nn
\end{align}
Since $rz^2\Rda$ is invertible, $T=\vN$ satisfies the
defining relations
\rf[plus]. Consequently, $\gdul$ is a quotient of the algebra $\gdsl$.
Putting both relations
together will complete the proof of (i).
\\
(ii) Throughout we sum over repeated
indices.
There are two $\adr$-invariant elements $W_\pm$ in
$\RR\cap\langle u^i_j,u^i_ju^n_l:\,i,\,j,\,n,\,l=1,\dots,N\rangle$, see
\cite{SS1}, formula (3) before Theorem\,2.1:
\begin{equation}
W_\pm:=\overline{V}_\pm-\mu_{\tau,k}^\pm\overline{U},\,\quad
 V_\pm=q^{-2j-2i}(P_\pm)^{ji}_{mn}u^m_ju^n_i,\quad U=\tsum
q^{-2i}u^i_i,
\end{equation}
where $\abar=a-\ve(a)$. The  numbers $\mu^\pm_{\tau,k}$ are given
in \cite{SS1}. Recall that $q_k$, $k=1,\dots,N$,  is an $N$th root of $q^2$.
Applying $\SS$ to $W_\pm$ one gets
\begin{equation}
\SS(W_\pm)=f_2^\pm\om_2 + f^\pm_{11}\om_1\ota\om_1,
\end{equation}
where $\om_2=\sum q^{-2i}\om^i_n\ota\om^n_i$. The constants can be
determined to be
\begin{align}
f^\pm_2&=q_k^{-2\tau}q^{\pm2(1+\tau)}(N\pm2)_q-\mu_{\tau,k}^\pm,
\\
f^\pm_{11}&=q^{3\pm1}(2)_q(N)_q^{-1}(N\pm 1)_q^{-1}(\mu^\pm_{\tau,k})^2-
q_k^{-2\tau}q^{\pm2(1+\tau)}(N)_q^{-1}(N\pm2)_q.
\end{align}
It takes several pages to prove the following fact: If $q$ is not a root
of unity and not a root of the algebraic equation
$q^{2\tau}p^N=(p+\rrr_\tau(N)_q)^N$, where $p=(N)_q^2+\rrr_\tau\QP^2(N)_q+
\rrr_\tau^2\QP^2$, then $f^+_2f^-_{11}-f^+_{11}f^-_2\ne0$.
In particular, if $q$ is  transcendental  the above determinant is
nonzero and there exist numbers $g^\pm_n$ with
$\om_1\ota\om_1=\SS(g^+_{11}W_++g_{11}^-W_-)$ and
 $\om_2=\SS(g^+_2W_++g_2^-W_-)$.
Hence, $\om_1\ota\om_1\in\SS(\RR)$ and
finally by \rf[zero], $\nu_1\ota\nu_1\in\SS(\RR)$.
\\
(iii)  The above proof fails in case $N=2$ since $\overline {V}_-=0$.
One can prove: If $\SS(r)$ is bi-invariant then one can find
an $\adr$-invariant element $r'\in\RR$ with $\SS(r)=\SS(r')$.
\\
For $\AA=SL_q(2)$,  $q$  transcendental, the space $\AAi$ of
$\adr$-invariant elements  is given by $\AAi=\langle U^k:k\in\NO\rangle.$
Set $\RRi=\RR\cap\AAi$. One can check that $\SS(\RRi)=\langle
2\om_1\ota\om_1+\tttt \om_2\rangle$,
where $\tttt=\pm q^{-5}(q\mp 1)(q^3\mp1)$ for the $4D_\pm$ calculus.
By Schur's lemma and irreducibility of $\adr$ on the subspaces
(1.24) and (1.25)
in \cite{Wo2} and since
\mbox{$\SS(b^2)=(q+\qm)\om(b)\ota\om(b)\ne0$} and
$\SS\bigl(\bigl(q^{-2}a+q^{-4}d\mp(\qm+q^{-5})\bigr)b\bigr)\ne0$,
$\SS$ is an isomorphism of $\adr$ and $\dr$ on the above  subspaces. Hence,
$\dim \SS(\RR)=9$.
Since $\SS(\RR)\subseteq\ker(I-\sig)$ the proof of (iii) is complete.
\end{Proof}

\end{document}